\documentclass[a4paper,10pt]{amsart}
\usepackage[english]{babel} 
\usepackage[leqno]{amsmath}
\usepackage{amssymb,amsthm}
\usepackage{amscd}
\usepackage{enumerate}

\newtheorem{thm}{Theorem}[section]
\newtheorem{lemma}[thm]{Lemma}

\newtheorem{prop}[thm]{Proposition}

\newtheorem{cor}[thm]{Corollary}

\theoremstyle{remark}
\newtheorem{remark}[thm]{Remark}
\theoremstyle{definition}
\newtheorem{defi}[thm]{Definition}

\newtheorem{example}[thm]{Example}

\newcommand{\la}{\longrightarrow}
\newcommand{\ha}{\hookrightarrow}
\newcommand{\da}{\dashrightarrow}

\newcommand{\ov}{\overline}

\newcommand{\Spec}{\operatorname{Spec}}

\newcommand{\im}{\operatorname{Im}}
\newcommand{\Part}{\operatorname{Part}}
\newcommand{\supp}{Supp}
\newcommand{\mdeg}{\operatorname{{\underline{de}g}}}

\def\X{\mathcal X}

\def\L{\mathcal L}

\def\O{\mathcal O}

\def\NN{\operatorname{N}}

\def\dcg{\Delta _X}

\def\int{M_X}

\def\md{\underline{d}}
\def\mt{\underline{t}}
\def\mdu{\underline{u}}

\newcommand{\Z}{\mathbb{Z}}

\newcommand{\mgbar}{\ov{M}_g}

\def\AN{\NN(\alpha^d_K)}

\newcommand{\uni}[1]{E_{#1}}

\newcommand{\dX}{\dot{\X}}
\newcommand{\dfib}{\dot{X}}

\newcommand{\partd}{{\operatorname{Part}}(d,\gamma)}

\newcommand{\DXO}{\mathcal{D}(X)}
\newcommand{\DX}{\overline{\mathcal{D}(X)}}
\def\nXr{N^{d, \rep}_X}

\def\nf{N^d_f}
\def\nfr{N^{d, \rep}_f}
\newcommand{\nq}{q_f}

\newcommand{\Pic}{\operatorname{Pic}}
\def\Picgen{\Pic_{K}}
\def\Picgend{\Pic_{K}^d}

\newcommand{\picf}[1]{\Pic_f^{#1}}

\newcommand{\picX}[1]{\Pic^{#1}X}
\def\JXbar{\overline{J_X^d}}
\def\JX{J_X^d}
\def\Jfbar{\overline{J_f^d}}
\def\Jf{J_f^d}

\newcommand{\tw}{\operatorname{Tw}}
\newcommand{\twf}{\tw_f}

\newcommand{\XS}{X_S^{\nu}}

\newcommand{\gen}{\X_K}

\newcommand{\rep}{\operatorname{r}}
\newcommand{\ar}{\alpha^{d,\rep}_f}
\newcommand{\ark}{\alpha^{d,\rep}_{f,X}}

\newcommand{\degt}{\Lambda_X^0}
\newcommand{\DXt}{\DX^0}
\newcommand{\ce}{{\epsilon}}
\newcommand{\pg}{\Part(d,\gamma)}

\newcommand{\singsep}{X_{\text{sing}}^{\text{sep}}}
\newcommand{\sing}{X_{\text{sing}}}

\begin{document}
\begin{center}
{\bf\large
Naturality of Abel maps} \\

%\bigskip

 Lucia  Caporaso

%\bigskip
%\today
\end{center}

\noindent
{\it Abstract.} We give a combinatorial characterization of nodal curves
admitting a natural (i.e. compatible with and independent of
specialization) $d$-th Abel map for any
$d\geq 1$.

%\tableofcontents

\

Let $X$ be a  smooth  projective curve and $d$ a positive integer; the classical $d$-th Abel
map of $X$,\  
$\alpha_X^d:X^d\la\Pic^dX$,
associates to $(p_1,\ldots,p_d)\in X^d$
the class of the line bundle  $\O_X(p_1+\ldots+p_d)$ in $ \Pic^dX$.
Such a morphism
 has  good functorial properties,
it is compatible with specialization and base change.

Now let $X$ be a  singular nodal curve occurring as the limit of a family of smooth curves.
We ask whether there is a   notion of   $d$-th Abel map for $X$ 
which is   limit of the Abel maps of the smooth curves of the family, and
which is {\it natural}, i.e. independent of the choice of the family.

It is known that, although a nodal curve $X$  is endowed with
a generalized Jacobian and a  Picard scheme which are both natural
(i.e. they are the limit of the Jacobians and, respectively, of the Picard schemes
of the fibers of every family of  curves specializing to $X$), there are
interesting
 degeneration problems
about line bundles and linear series where naturality is a subtle issue.
Abel maps is the case which we investigate here.

The main result of this paper, Theorem~\ref{main},
characterizes in purely combinatorial terms nodal curves that   possess a natural $d$-th
Abel map. 

A consequence of our result is that, if we consider 
stable curves of genus $g\geq 2$, then the locus in  $\mgbar$ of curves that fail to admit
a natural $d$-th Abel map, for a fixed $d\geq 2$, has codimension $2$.
So, naturality of  Abel maps is not  to be taken for granted, unless $X$
is irreducible
or of compact type, in which case it is not difficult to see that  natural Abel maps exist
for all
$d$.

What is a good notion of  Abel maps for singular curves?
The same definition as for the smooth case
 behaves badly under specialization; moreover,
it obviously does not make sense if some of the $p_i$ are singular points of $X$.
This last problem will not be an  issue here:
 we shall only study  noncomplete Abel maps, which is enough
for our scopes. 

So, for us a $d$-th Abel map is a rational map $\beta:X^d\da \Pic^dX$
arising  as the limit  of the   Abel maps  of smooth curves 
specializing to $X$, for some family.
A further requirement is added to ensure separation;
 see  \ref{abdef}.

Thus, the target space of our $d$-th Abel map is the  Picard scheme, not 
 any particular 
compactification of it. 
 Our definition  and  results  should be sufficently general  to apply to
various
   compactified Picard schemes existing in the literature (see section~\ref{compare}). 

The construction of complete Abel maps for singular curves was 
carried out
 by A. Altman and S. Kleiman for irreducible and reduced  curves
in \cite{AK}; see also \cite{EGK} for further results.
 Not  much is known for reducible curves. Recently, in \cite{CE},
 degree-$1$ Abel maps of stable curves  are defined,  compactified,
 and shown to be  natural.
  For higher $d$
  the completion problem is   open in general, see \cite{coelho} for some
progress in case $d=2$.  
The main result of the present paper  indicates that a safe way to approach it
is to work with a fixed  one-parameter smoothing of the given curve
$X$ (as in \cite{CE} and \cite{coelho}), or to restrict to   natural Abel
maps.
     
 Among our techniques, the main one is the    
 use N\'eron models of Jacobians (as constructed by M. Raynaud in \cite{raynaud}); this allows us to
obtain a concrete description of our  axiomatically defined Abel maps.
Then we combine a result of E. Esteves and N. Medeiros
 about  deformation  of line bundles and enriched
structures (in \cite{EM}) with
a detailed combinatorial analysis.

I wish to thank Simone Busonero and Eduardo Esteves for  precious suggestions
and remarks.

\section{Statement of the main result}
\label{}
In this section we state the main theorem
(\ref{main}) after a few preliminaries.

\subsubsection{Conventions}
\label{not}
We work over   an algebraically closed field $k$.
$X$  always denotes a connected, reduced,  projective curve defined over $k$
having at most nodes as singularities, and 
 $C_1,\ldots, C_{\gamma}$ its 
 irreducible components. 

By a (one-parameter) {\it regular smoothing of $X$} we  mean
a proper morphism $f:\X \to  B=\Spec R$,
with $R$   a discrete valuation ring
having   residue field $k$ and quotient field $K$,
such that
   $X$ is the closed fiber of $f$, the total space $\X$ is nonsingular, and such that
the  generic fiber of $f$, denoted
$\X_K$, is a smooth projective curve over $K$.

Let $X=\cup_1^{\gamma}C_i$ be a curve as above and $L\in \Pic X$ a line bundle on it. 
Then $L$ has a {\it degree}, $d=\deg L \in \Z$, and a 
{\it multidegree}, $\md=\mdeg L \in \Z^\gamma$ defined as
$\mdeg L:=(\deg_{C_1}L,\ldots,\deg_{C_{\gamma}}L )$.
Denoting $\Pic^dX$, respectively $\picX{\md}$, the locus in $\Pic X$ of line bundles of degree
$d$, respectively of multidegree $\md$, we have
$$
\Pic X=\coprod _{d\in \Z}\Pic^dX
\hskip.7in \Pic^d X=\coprod _{\md\in \Z^{\gamma}:|\md|=d}\picX{\md},
$$
where $|\md|:=\sum_1^{\gamma} d_i$ for every $\md=(d_1,\ldots,d_{\gamma}) \in \Z^\gamma$.
$\Pic^{\md}X$ is isomorphic to the generalized Jacobian of $X$.

If $Z\subset X$ is a (complete) subcurve,
we denote  by $\md _Z$ the multidegree of $L_{|Z}$ ($L$ restricted to $Z$)
and by $|\md _Z|$ the total degree of $L_{|Z}$.

Set
$
\DXO: =\{ D=\sum n_i C_i,\  n_i\in \Z\}
$
 the free abelian group  generated  by $C_1,\ldots, C_{\gamma}$.
If $Z\subseteq X$ is a   subcurve of $X$, so that $Z=\cup_{C_i\subset Z}C_i$, we
shall, as usual, abuse notation by denoting again $Z=\sum_{C_i\subset Z}C_i\in \DXO$.

For any $f:\X \la \Spec R$
be a regular smoothing of $X$
we have a
symmetric bilinear product  $(\ \cdot\ ):\DXO\times \DXO\to \Z$, 
often called  the ``intersection pairing",
which is the same for every $f$
(as long as the total
space $\X$ is regular, which we always assume).  
Recall  that
$(X\cdot D)=0$ for  all $D\in \DXO$.

For any subcurve $Z\subset X$ 
we set  $Z':=\overline{X\smallsetminus Z}$
and 
$
k_Z:=(Z\cdot Z')=-Z^2
$.

Consider the quotient $\DX$ of $\DXO$ by the subgroup generated by $X$:
$
\DX :=\frac{\DXO}{\Z \cdot X}.
%=\frac{\{\sum n_iC_i \  n_i\in \Z\}}{\Z \cdot X}.
$
The intersection pairing descends to 
$\DX$.

\subsubsection{Twisters}
\label{twnot}
For a fixed  regular smoothing $f:\X \to B$  of $X$,  
set
$$
\twf X:=
\frac{\{\O_\X(D)_{|X}, \   \forall D\in \DXO\}}{\cong}
\subset \Pic^0X.
$$
The union as $f$ varies among the regular smoothings of $X$ is denoted  by
\begin{equation}
\label{twist}
\tw X := \bigcup _{f \text{ reg. sm.}}\twf X\subset \Pic^0X
\end{equation}
Elements of $\tw X$  are special line bundles called {\it twisters}.
%(see \cite{cner} 3.2)
Any  $T\in \tw_fX$ determines a  $D\in \DXO$ such that
$T\cong\O_\X(D) _{|X}$ only up to adding a multiple of the central fiber $X$ of $f$; 
in fact
$\O_\X(D) _{|X}\cong \O_\X(D+nX) _{|X}$ for any $n\in \Z$.
Thus  we have a surjective map
$
\tw X \la \DX
$
associating to any $T\in \tw X$  the class of a  $D\in \DXO$ such that
$T=\O_\X(D) _{|X}$.
We shall denote such a class $\supp (T)$ and call it the {\it  support of $T$}.

Let $D\in \DXO$ and $T\in \tw X$    such that $\supp (T)=\ov{D}$,
then 
$
\mdeg T=((D\cdot C_1), \ldots, (D\cdot C_{\gamma}))
$,
independently of the  representative $D$ for $\ov{D}$.
Furthermore, this shows that $\mdeg T$ does not depend on the regular smoothing
defining $T$. In other words, the multidegree of a twister only depends on its support,
so that
 we can 
unambiguously write
\begin{equation}
\label{degsum}
\mdeg D:= \mdeg \O_\X(D) _{|X}=((D\cdot C_1), \ldots, (D\cdot C_{\gamma}));
\end{equation}
for a  subcurve $Z\subset X$ we set 
$
\mdeg _ZD :=\mdeg \O_\X(D) _{|Z} 
$.
Summarizing, we have 
$$
\mdeg:\tw X \stackrel{\supp}{\  \la\  } \DX \stackrel{\mdeg}{\  \la\  } \Z^{\gamma}
$$
where  the second arrow  is a group homomorphism.
We denote by $\Lambda _X \subset \Z^{\gamma}$ the image of $\mdeg$ above, 
i.e. the group of multidegrees of all twisters:
\begin{equation}
\label{Lambda}
\Lambda _X:=\{\md\in \Z^{\gamma}: \exists T\in \tw X: \mdeg T=\md\}.
\end{equation}
More details about this set up will be in section~\ref{tw}.

Interpreting $\Z^{\gamma}$ as the set of all multidegrees on $X$, we get an equivalence
relation on it, induced by $X$:

\noindent
{\it Let $\md, \md'\in \Z^{\gamma}$, define  $\md$ equivalent to $\md'$
 (in symbols $\md\equiv \md '$ )
iff $\md-\md '\in \Lambda_X$}

\noindent
where $\Lambda_X$ is defined in (\ref{Lambda}).
Introduce the set of multidegree classes of total degree equal to a fixed $d$ 
\begin{equation}
\label{dcg}
\dcg^d:=\frac{\{\md \in \Z^{\gamma}: |\md |=d\}}{\equiv}.
\end{equation}
 $\dcg^d$ is well known 
to be a finite set
whose cardinality does not depend on $d$
(\cite{raynaud} 8.1.2, see also \cite{cner} 3.7 for an overview). 
As we shall see,   $\dcg^d$ is useful to control the
non-separatedness of the Picard scheme.

For $\md \in \Z^{\gamma}$ we shall denote $[\md] \in \dcg^d$ its class.

\subsubsection{}
\label{pic}
Let $f:\X \to B$ be a regular smoothing of $X$.
Then there exists a
Picard scheme relative to $f$,
denoted here
$\Pic _f$ (an alternative    notation is $\Pic _f = \Pic _{\X/B}$, which is not used in this paper).
$\Pic_f$ is a scheme over $B$ whose  generic fiber  is
$\Pic_{\gen}$,
 the Picard scheme of the curve $\gen$.
$\Pic _f$
has a basic moduli property which we need to recall. Let $S$ be a 
$B$-scheme, then for any line bundle $\L$ on $S\times _B\X$ there exists a unique $B$-morphism
\begin{equation}
\label{modL}
\mu_\L : S \la \Pic_f\  \  \  s\mapsto \L_{|k(s) \times _B\X}
\end{equation}
which we refer to as the {\it moduli map of } $\L$. More details are in \cite{GIT} 05
(d).

$\Pic_f$ decomposes  into its connected
components, 
$\Pic_f=\coprod_{d\in \Z} \picf{d}$, 
where $\picf{d}\to B$ parametrizes line bundles of relative degree
 $d$. The generic fiber of $\picf{d}\to B$ is denoted
$\Picgend:=\Pic_{\gen}^d$, 
an integral projective variety over $K$; 
the closed fiber is $\Pic^dX$, 
a  reducible scheme if $X$ is reducible (see \ref{not}).

$\Pic_f$ is  smooth    over $B$
but    not separated if $X$ is reducible.
The essential reason  is  the existence
of twisters (see \ref{twnot}), which  must all be identified in
any separated completion of
$\Picgen$ over $B$.

Using the moduli property of $\Pic_f$ one can construct Abel maps 
for smooth curves
over any base scheme, see \ref{abel2}. For example
let $X$ be a smooth curve over $k$. As we already said, the $d$-th Abel map of
$X$ is the map
$\alpha_X^d:X^d\to \Pic^dX$ such that
$\alpha_X^d(p_1,\ldots, p_d)=\O_X(p_1+\ldots + p_d)$.

We  shall denote $\alpha ^d_K$  the $d$-the Abel map of the generic fiber of $f$:
\begin{equation}
\label{abgen}
\alpha ^d_K =\alpha ^d_{\gen}: \gen ^d \la \Picgend
\end{equation}
\begin{example}
\label{example}
Let $X=C_1\cup C_2$ be a  curve having two smooth components meeting in only one node $r$.
Let us examine the naive definition for the Abel map in degree $1$, 
  copying the  smooth case. We  get a map
regular away from the node $r$:
\begin{equation}
\label{naive}
\begin{array}{lccr}
\dfib &\stackrel{\alpha}{\la}&
\Pic^1X&=\coprod_{d_1+d_2=1}(\Pic^{d_1}C_1\times \Pic^{d_2}C_2)\\  p &\mapsto
& \O_X(p)& 
\end{array}
\end{equation}
where $\dfib:=X\smallsetminus \{r\}$.
Let us illustrate some pathologies of  definition (\ref{naive}).
 The two  components of $\dfib$ 
are mapped to two different connected components of $\Pic^1X$,
namely  $\alpha (C_1)\subset \Pic^1C_1\times \Pic^0C_2$
while $\alpha (C_2)\subset \Pic^0C_1\times \Pic^1C_2$.
All the connected components $\Pic^{d_1}C_1\times \Pic^{d_2}C_2$
of $\Pic^1X$ are obviously projective, hence $\alpha$ cannot possibly be extended to
a regular map from the whole of $X$ to $\Pic^1X$. 

A second problem is  the fact that
$(1,0)\equiv (0,1)$; indeed $\#\dcg^1 =1$  (see
(\ref{dcg})).   
Now   if $\X\to B$ is a regular smoothing of $X$,  
any separated model of $\Pic^1_{\gen}$ over $B$ cannot contain components of $\Pic X$
corresponding to equivalent multidegrees (see below).
So
the target  of the  naive Abel map of the family fails to be separated.
\end{example}

A more satisfactory definition turns out to be \ref{abdef} below; first
some notation.
Let $f:\X \to B$ be a regular smoothing of  $X$, then
 $f_d:\X^d_B\to B$ denotes the $d$-th fibered power of $\X$ over $B$.
The open subset of $\X^d_B$ where  $f_d$ is smooth is denoted by
$
\dX^d_B:=\X^d_B\smallsetminus sing(f_d).
$
Similarly, $\dfib ^d := \{(p_1,\ldots,p_d): p_i\in X\smallsetminus \sing\}
$ denotes the closed fiber of $\dX^d_B\to B$,
 where $\sing$
is the set of   singular points  of  $X$.

\begin{defi}
\label{abdef}
%\begin{enumerate}[(1)]
%\item
\label{}
A {\it $d$-th Abel map} for the curve $X$ is a regular map
$
\beta:\dfib^d\la \Pic^dX
$
satisfying the following requirements.
\begin{enumerate}[(a)]
\item
\label{spe} 
There exist a regular smoothing $f:\X\to B$ of $X$ and a map
$\beta_f:\dX ^d_B\la \picf{d}$ such that 
 the restriction of $\beta_f$ to the generic fiber is the $d$-th Abel map of
$\gen$ (i.e. $(\beta_f)_{K}=\alpha^d_K$ of (\ref{abgen})),
and the restriction of $\beta_f$ to the closed fiber is equal to
$\beta$ (i.e. $(\beta_f)_k=\beta).$
\item
\label{sep}
If $\md,\md'\in\Z^{\gamma}$ are such that $\im\beta\cap \Pic ^{\md}X\neq \emptyset$ 
and $\im\beta\cap \Pic ^{\md'}X\neq \emptyset$, then $\md \not\equiv \md'$.
\end{enumerate}
%\item
\label{}
A $d$-th Abel map $\beta$ is called {\it natural}
if it is independent on the choice of $f$.
More precisely, if for every regular smoothing $f$ of $X$ there
exists a  (necessarily unique) $\beta_f$
as in (\ref{spe}), extending $\beta$.
%\end{enumerate}
\end{defi}
Condition~(\ref{sep})  will ensure that the image of an Abel map
is contained  in some 
separated model of $\Picgend$ over $B$ (cf. Example~\ref{example}).

A  concrete description of Abel
maps  will be    given by Proposition~\ref{key}.
In section~\ref{compare} we shall relate our definition to others in the
literature.

Although (see \ref{exist})
every curve   does have  Abel maps,  not all curves
admit a natural one. Our main result, Theorem~\ref{main} characterizes 
in purely combinatorial terms those curves   admitting a natural $d$-th Abel map. 
Before stating it we   define the crucial combinatorial character.
Denote 
$\singsep\subset \sing$
  the set of   {\it separating  nodes} of $X$,
 i.e.  $\singsep$ is the set of
 nodes $p$ of $X$  such that  $X\smallsetminus p$ is
disconnected. 

In the next definition, we use the notation of \ref{not} and adopt the convention
that if $S$ is an empty set of integers, then   \  
$\inf \{n\in S\}=+\infty$.

\begin{defi}
\label{kconn}
The {\it essential connectivity}  $\ce(X)$ of   $X$  is 
\begin{equation}
\label{ec}
\ce(X):=\inf \{k_Z, \  \forall Z:\  \emptyset \subsetneq Z \subsetneq X \text{ and } Z\cap
Z'\cap\singsep=\emptyset
\}.
\end{equation}
Equivalently
\begin{equation}
\label{ec2}
\ce(X):=\inf \{k_Z, \  \forall Z:\  \emptyset \subsetneq Z \subsetneq X \text{ and } Z\cap
Z'\not\subset\singsep
\}.
\end{equation}
\end{defi}

\begin{remark}
\label{}
So, if $X$ is irreducible 
or  of compact type  then
$\ce(X)=+\infty$.

An elementary arguments yields that   to compute 
 $\ce(X) $ it suffices to consider connected subcurves  and that
 (\ref{ec}) and (\ref{ec2}) are equivalent; we omit it since  we shall 
only use version (\ref{ec2}).
An example: if $X=C_1\cup C_2$ then either  $\#(C_1\cap C_2)= 1$ and
$\ce(X)=+\infty$, or
$\#(C_1\cap C_2)\geq 2$
and $\ce(X)=\#(C_1\cap C_2)$.
\end{remark}

\begin{thm}
\label{main} 
Let $X$ be a nodal curve and  $d$  a positive integer. Then $X$ admits a natural $d$-th Abel map
 if and only if $\ce (X)>d$.
\end{thm}
The proof of the theorem will be in Section~\ref{proof}.

\section{Abel maps via N\'eron models}
The goal of this section is to obtain a complete description of Abel maps, which will
be done in \ref{key}.
We shall use N\'eron models  in the same spirit of \cite{cner} (and of \cite{edix} section 9). We
 refer to \cite{BLR} chapter 9
or to \cite{artin} section 1 for details.

Let $f:\X \to B$ be a regular smoothing of $X$; as we said in \ref{pic},
if $X$ is reducible $\picf{d}$ is not separated over $B$.
To correct this  we
introduce  
the N\'eron model 
$\NN(\Picgend)\to B$ 
of  its generic fiber $\Picgend$. This is a smooth, separated scheme over $B$  uniquely determined by
a universal property, the N\'eron mapping property
(\cite{artin} 1.1 or \cite{BLR} 1.2.1).
 We shall denote
$$
N^d_f:=\NN(\Picgend).
$$
Recall that there exists a  unique surjective $B$-morphism
\begin{equation}
\label{quotmap}
\nq:\picf{d}\la N^d_f
\end{equation}
which is the identity on the generic fiber
($\nq$ is  called ``Ner" in \cite{artin} diagram 1.21).

%\subsubsection{}
To describe  the map (\ref{quotmap}),
consider $\picf{\md}\subset \picf{d}$ the moduli scheme for degree-$d$
line bundles
having 
multidegree $\md$ on the closed fiber $X$ of $f$.
Then
$$
\picf{d}=\frac{\coprod _{|\md|=d }\picf{\md}}{\sim_K}
$$
where ``$\sim_K$"  denotes the natural gluing of the schemes $\picf{\md}$ along their
generic fiber (which is the same for all $d$: $\Picgend$). 
Similarly  we have
\begin{equation}
\label{nerdesc}
\NN^d_f\cong\frac{\coprod _{\delta \in \dcg^d}\picf{\delta}}{\sim_K}
\end{equation}
where 
$
\picf{\delta}:=\picf{\md}
$ for any  representative $\md$ of $\delta$.
Such a definition does  not  depend on the choice of $\md$;
in fact
 for every pair of equivalent multidegrees  $\md, \md'$ there exists a unique
isomorphism between $\picf{\md '}$ and $\picf{\md}$, determined by
the  unique $T\in \tw_fX$ such that $\mdeg T=\md - \md '$
(see \ref{twlm} (\ref{twuni})). 

The  map
$\nq$  restricted to 
$\picf{\md}\subset \picf{d}$ is the unique isomorphism 
$
\picf{\md}\cong \picf{[\md]}.
$

We denote $N_X^d$ the closed fiber  of $\NN^d_f\to B$, which, if $X$ is reducible,  is  
the  disjoint union of $\#\dcg^d$  copies of the generalized  Jacobian of $X$.

Although $\nf$  does not have good functorial properties,
  it has a geometric interpretation.
Let 
$L,L'\in \Pic X$ and call them 
{\it $f$-twist equivalent} iff $L'\otimes L^{-1}\in \tw_fX$.
Then
 the fibers of 
$
(\nq)_k:\picX{d}\la N^d_X
$
are the  classes of $f$-twist equivalent 
line bundles.

\subsubsection{Abel maps and N\'eron models}
\label{abel2}
To  study Abel maps of singular curves, we 
use an  approach analogous to \cite{GIT} section 6, just like in \cite{CE}.
With the notation introduced in \ref{pic}, let $f:\X \to B$ be a regular smoothing
of  $X$. Consider the base change of $f$ to  $\dX ^d_B\to B$, namely
$$\pi:\dX ^d_B\times _B\X\la \dX ^d_B$$
so that $\pi$ is the first projection.
Define the ``universal effective (Cartier) divisor"
$\uni{d}$ on $\dX ^d_B\times _B\X$ as the sum of the $d$ natural sections 
$\sigma_1,\ldots, \sigma_d$ of
$\pi$, that is, the sections 
$
\sigma _i (p_1,\ldots,p_d)=((p_1,\ldots,p_d),p_i).
$
Consider the moduli map of $\uni{d}$, 
$
\mu_{\uni{d}}:\dX ^d_B\to \picf{d}
$
as defined in (\ref{modL}). By definition,
its restriction   to the generic fiber 
is its  $d$-th Abel map,  
$
\alpha ^d_K : \gen ^d \to \Picgend,
$ introduced in (\ref{abgen}).
Consider now the composition
\begin{equation}
\label{Nmap}
\AN:\dX ^d_B \stackrel{\mu_{\uni{d}}}{\la}\picf{d}\stackrel{\nq}{\la} N^d_f
\end{equation}
where the 
 notation $\AN$ is motivated as follows.
To
the  map $\alpha ^d_K $ 
  one can apply  the N\'eron mapping property:
since $\dX ^d_B$ is smooth over $B$ and has $\gen^d$ as generic fiber, there exists a
unique extension of 
$\alpha ^d_K $ to a morphism $\dX ^d_B\to \nf$; such an extension, 
denoted $\AN$, is unique and hence  coincides
with the  map (\ref{Nmap}).

Finally, notice that if $D$  is a Cartier divisor  supported on the closed fiber of $\pi$
(i.e. on $\dfib^d\times X$), the same construction,   replacing
$\uni{d}$ with $\uni{d}+D$,  gives again $\AN$
(in other words: $\AN= \nq\circ\mu_{\uni{d}+D}=\nq\circ\mu_{\uni{d}}$).

To  use the morphism (\ref{Nmap})
we shall     fix   a choice of representatives for the multidegree classes,
given by a map
$$
\rep : \dcg^d \la \Z ^{\gamma}\  \text{ such that }\       
|\rep (\delta ) |=d \   
  \text{ and } \  [\rep(\delta)]=\delta,\     
\forall \delta\in \dcg^d.
$$
 Then there
exists a unique isomorphism 
\begin{equation}
\label{rep}
\iota_{\rep}:\nf \stackrel{\cong}{\la} 
\frac{ \coprod_{\delta \in \dcg^d} \Pic^{\rep(\delta)}_f}{\sim_K}=:\nfr
\end{equation}
So $\nfr$ defined above is a subset of $\picf{d}$. 
Restricting $\nq:\Pic_f^d\to\nf$ to it we have
 $id_{\nfr}=\iota_{\rep} \circ(\nq)_{|\nfr}$.
The closed fiber of $\nfr$ is denoted
$$
\nXr:=(\nfr)_k=\coprod_{\delta \in \dcg^d} \Pic^{\rep(\delta)}X\subset \Pic^dX
$$
and obviously does not depend on $f$.
Now compose the map $\AN$ of (\ref{Nmap}) with  $\iota_{\rep}$ and call the composition
  $\ar$
\begin{equation}
\label{ar}
\ar: \dX^d \stackrel{\AN}{\la}\nf \stackrel{\iota_{\rep}}{\la}\nfr .
\end{equation}
Restricting  to the closed fiber we get a $d$-th Abel map
$
\ark:\dfib^d\la \nXr
$.

To complete the picture,
 focus on the  multidegrees that are attained by the
divisor
$E_d$ restricted to the fibers of $\dX^d_B\times_B \X\to \dX^d_B$.
We shall name them {\it partitional multidegrees} and denote their set $\pg$:
$$
\partd:=\{\md \in \Z^{\gamma}: |\md |=d \text{ and }  \  d_i\geq 0 \   \forall i=1,\ldots \gamma \}
$$
Then the  map $\ar$ factors
\begin{equation}
\label{arfact}
\ar:\dX ^d \stackrel{\mu_{E_d}}{\  \la\  } \frac{ \coprod_{\md \in \pg}
\Pic_f^{\md}}{\sim_K}
\stackrel{q_f}{\la}
\nf\stackrel{\iota_{\rep}}{\la}\nfr
\end{equation}
where, abusing notation,
$q_f$ is the restriction of  $q_f :\Pic^d_f\la \nf$.
Therefore the image of $\ark$   is entirely contained in the union of the
 components  of $\nXr$ that correspond to those
multidegree classes  containing some
partitional representative. 
%(as explicitely stated in  (\ref{ark}) below).
The conclusion of the preceeding discussion is stated in part (\ref{key1})
of \ref{key} below,  where we also classify   Abel maps as defined in \ref{abdef}.
\begin{prop}
\label{key} 
\begin{enumerate}[(i)]
\item
\label{key1}
Let $f:\X \to B$ be a 
regular smoothing of  $X$, $d$ a positive integer and $\rep:\dcg^d \la \Z^{\gamma}$
a choice of representatives. Then there exists a unique morphism
$
\ar: \dX ^d_B \la \nfr
$
whose restriction to the the generic fiber is the  $d$-th Abel map of $\gen$.
The restriction of $\ar$ to the closed fiber is the $d$-th Abel map
\begin{equation}
\label{ark}
\ark:\dfib^d\la \coprod_{\md \in \pg} \Pic^{\rep([\md ])}X\subset \nXr .
\end{equation}
\item
\label{}
Conversely, every $d$-th Abel map of $X$  equals  $\ark$ for some  $f$ and $\rep$.
\end{enumerate}
\end{prop}
\begin{proof} The first part has  been proved before the statement.  
For the
second, let $\beta:\dfib^d\to \Pic^dX$ be an Abel map as defined in \ref{abdef}.
Then there exists a regular smoothing
$f:\X \to B$ and a morphism $\beta_f:\dX ^d_B\to \picf{d}$ extending $\beta$ and
restricting to the $d$-th Abel map $\alpha_K^d$ of $\gen$ on the generic fiber.

Consider the composition 
$
\nq\circ\beta_f:\dX^d_B\la \picf{d}\la \nf
$.
We claim that
\begin{equation}
\label{qb}
\nq\circ\beta_f=\AN
\end{equation}
where $\AN$ is defined in (\ref{Nmap}). Indeed the two maps $\nq\circ\beta_f$
and $\AN$
coincide on the generic fiber $\gen ^d$ and hence, by the
unicity part in the N\'eron mapping property,
they are equal.

Now define   $\rep:\dcg^d\to \Z^{\gamma}$ as follows;
pick $\delta \in \dcg^d$. 
If there exists a representative $\md^{\delta}$ for $\delta$
such that
$\im \beta \cap \picX{\md^{\delta}}\neq \emptyset$, then such a $\md^{\delta}$
is unique by condition (\ref{sep}) of definition~\ref{abdef}; hence we can define
$\rep({\delta})=\md^{\delta}$.
If instead there is no such   representative for $\delta$, we define $\rep$ however we
like.

By construction
$$
\im \beta_f\subset \frac{\coprod_{\delta \in \dcg^d}
\Pic^{\rep(\delta)}_f}{\sim_K}=\nfr
$$
Now consider the isomorphism $\iota_{\rep}$ (see (\ref{rep}))
$
\iota_{\rep}:\nf \stackrel{\cong}{\la} 
\nfr\subset \picf{d}.
$

Recall that
$\iota_{\rep} \circ(\nq)  _{|\nfr}=id_{\nfr}$.
Therefore
$$
\beta_f=id_{\nfr}\circ\beta_f=\iota_{\rep} \circ(\nq) 
_{|\nfr}\circ\beta_f=
\iota_{\rep} \circ\AN 
$$
using (\ref{qb}) for the last equality.
Now $\iota_{\rep} \circ\AN=\ar$  by definition (see (\ref{ar})),
therefore  $\beta =(\beta_f)_k=\ark$ and  we are done.
\end{proof}

\begin{remark}
\label{exist}
We get that,  for all $d$, every curve $X$ has  Abel maps, infinitely many of them
if $X$ is reducible (at least one for
every  $\rep$).
%%!For example, the Abel maps considered in \cite{CE}
%correspond to choosing in  each class a ``balanced" representative, whenever it exists
%(\cite{CE} 3.2).
\end{remark}
As a consequence of \ref{key}, to prove Theorem~\ref{main} it suffices to study the maps $\ark$.
If $\rep$ is fixed, the domain and the target 
$\coprod_{\md \in \pg} \Pic^{\rep([\md ])}X$ of $\ark$ do not depend on $f$, whereas
the map itself does.
In fact $\ark$ factors through certain isomorphisms $\picX{\md}\to \picX{\rep([\md])}$;
each of these isomorphisms is  given by tensor product by the (unique)
$T\in \tw_fX$ of multidegree $\rep([\md])-\md$.
The crux of the matter is that $T$ may change as $f$ varies
(even though its multidegree and support do not change),
in which case the Abel map $\ark$ will not be natural.

\section{Spaces of twisters}
\label{tw}
The goal of this section is to characterize twisters that depend only on their support,
  hence only on their multidegree (see \ref{twlm}),
and not on the regular smoothing defining them (see \ref{indtwi}).
Recall  the set up   of \ref{twnot}; we shall need the following well known facts, see
 for example \cite{artin}  p.220 diagram 1.21. The notation is in \ref{twnot}. 
\begin{lemma}
\label{twlm}
\begin{enumerate}[(i)]
\item
\label{degtw} 
The map $\mdeg:\DX \la \Lambda _X$ is an isomorphism
\item
\label{twuni}
For any regular smoothing $f$ of $X$ and any $ T,T'\in \twf X$ we have 
$
T=T'\Leftrightarrow
\mdeg T=\mdeg T'$. In particular $\tw_fX\cong \DX$
\end{enumerate}
\end{lemma}

\begin{remark}
\label{Df}
For any  $D\in \DXO$ 
and   $f:\X \la B$ regular smoothing of $X$, there exists a unique 
$T\in \tw_fX$ such that $\supp(T)=\ov{D}$.
We will denote by $D_f$ such a $T$, that is
$
D_f\cong \O_\X(D) _{|X}
$.

Denote by $\mt:=\mdeg D$; the lemma above implies that  the set $\tw_fX$ contains a
unique 
 element of multidegree $\mt$.
By contrast the set
\begin{equation}
\label{twt}
\tw^{\mt}X:=\{T\in \tw X: \mdeg T = \mt\}=\{D_f, \  \forall f \text{ reg. sm. of } X\}
\end{equation}
may be quite large.
If  $\#\tw^{\mt}X=1$, i.e. if
$D_f=D_{f'}$ for all regular smoothings $f$, $f'$ of $X$, we will say that
{\it $D_f $ does not depend on the choice of $f$.}
\end{remark}

\subsubsection{Level curves}
\label{level}
Let  $D=\sum n_iC_i$, 
 we can write 
\begin{equation}
\label{Ddec}
D=\sum_{m\in \Z} mD_m 
\end{equation}
with  $D_m:=\cup _{n_i=m}C_i$  
(in particular  $D_0:=\sum_{n_i=0} C_i$)
so that the $D_m$ are possibly empty subcurves
of $X$ having no components in common, and  such that
$
\cup _{m\in \Z}D_m=X.
$ 
We call (\ref{Ddec}) the {\it level expression} of $D$ and the 
non-empty curves $D_m$
 the {\it level curves} of $D$; 
of course they are uniquely determined.

Notice also that for any $n\in \Z$ the level curves of $D$ and of $D+nX$  are the same,
hence it makes sense to  speak of level curves of a class $\ov{D}\in \DX$.
We can  also define level curves of a twister $T\in \tw X$
 as the level curves of its support, $\supp (T)\in \DX$.
Similarly, for any $\mt \in \Lambda_X$
the level curves of $\mt$ can be defined via the isomorphism $\DX\cong \Lambda _X$ of \ref{twlm}.
More precisely, we need the following
\begin{lemma}
\label{lvlm} 
Let $\mt \in \Lambda_X$ with $\mt\neq 0$. There exists a unique $D(\mt)\in \DXO$,
with
$\mdeg  D(\mt)=\mt$,  admitting an expression 
\begin{equation}
\label{Duni}
D(\mt)=\sum_1^{\ell (\mt)} m_hZ_h(\mt)
\end{equation}
where $\ell (\mt)$, $m_h$ and $Z_h(\mt)$ are uniquely determined by the following 
properties.
\begin{enumerate}[(a)]
\item
\label{a}
$\ell(\mt)\geq 1$ and  $m_h\in \Z$ with $0<m_1<\ldots<m_{\ell (\mt)}$; 
\item 
\label{b}
the $Z_h(\mt) $  are subcurves of $X$ having no components in common;
%(we  call them the {\em level curves of} $\mt$);
\item
\label{c}
the curve $Z_0(\mt):=\overline{X\smallsetminus \cup_1^{\ell(\mt)} Z_h(\mt)}$ 
 is not empty.
\end{enumerate}
For every subcurve $Y\subseteq Z_0(\mt)$ we have
\begin{equation}
\label{degZ}
|\mt_Y|=\deg_{Y}D(\mt) \geq -m_1(Y\cdot Z_0(\mt))\geq 0.
\end{equation}
In particular, if $Y=Z_0(\mt)$ we have $|\mt_Y|\geq m_1k_{Z_0(\mt)}>0$.

\end{lemma}
\begin{proof}
Pick any $D\in \DXO$ with $\mdeg D=\mt$ and consider the level expression
of $D$  in
(\ref{Ddec}). It has a finite number of nonzero summands, so let $m_0$ be the minimum integer
for which $D_{m_0}$ is not empty. Set $D(\mt ):=D-m_0X$. Then the level expression of $D(\mt)$
can be written as in (\ref{Duni}) and
 satisfies the conditions (\ref{a}), (\ref{b}), (\ref{c}) in the statement; note that $D_{m_0}=Z_0(\mt)$. 
 We need to check condition
(\ref{degZ}). So let $Y\subseteq Z_0(\mt)$, then
$$
|\mt_Y|=(Y\cdot D(\mt))=\sum_1^{\ell (\mt)} m_h(Y\cdot Z_h(\mt));
$$
now for every $h\geq 1$  we have $(Y\cdot Z_h(\mt))\geq 0$ 
(because $Y$ lies in the complementary curve
of
$Z_h(\mt)$). Therefore, since $m_1\leq m_h$ if $h\geq 1$, we get
$$
|\mt_Y|\geq m_1\sum_1^{\ell (\mt)} (Y\cdot Z_h(\mt))=m_1(Y\cdot \cup_1^{\ell(\mt)} Z_h(\mt))=
-m_1(Y\cdot Z_0(\mt))\geq 0
$$
indeed the $Z_h(\mt)$ have no common components and
 $\cup_1^{\ell(\mt)} Z_h(\mt)$ is a reduced proper subcurve of $X$
whose complement is, by definition, $Z_0(\mt)$. If $Y=Z_0(\mt)$ the last inequality is strict,
so we are done.
\end{proof}
Let $D\in \DXO$; to the  level expression $D=\sum_{m\in \Z}mD_m$ of
(\ref{Ddec}) we can naturally associate a  set of nodes $S(D)$ as follows. 
\begin{equation}
\label{S(D)}
S(D):=\cup_{m\neq m'}(D_m\cap D_{m'})\subset \sing
\end{equation}
Similarly, let $\mt\in \Lambda_X$ and consider
$D(\mt)=\sum_1^{\ell (\mt)} m_hZ_h(\mt)$ of \ref{lvlm}.
Denote
\begin{equation}
\label{S(t)}
S(\mt):=S(D(\mt))=\bigcup _{0\leq h< h'\leq \ell(\mt)}(Z_h(\mt)\cap Z_{h'}(\mt))\subset X_{sing}
\end{equation}
%Finally, for  $T\in \tw X$ we set $ S(T):=S(\mdeg (T)). $

\begin{remark}
\label{sumD}
For  $D, D'\in\DXO$
we have $S(D+D')\subset S(D)\cup S(D')$. 
%
%This is straightforward: it suffices to prove it for $D'=C$ an irreducible component of $X$, in
%which case $S(C)=C\cap C'$.
%We have
%$$S(D+D')=S(D+C)\cap \sing= S(D+C)\cap (C'_{\text sing} \cup (C\cap C')).$$
%Obviously $S(D+C)\cap C'_{\text sing}\subset S(D)$, whereas
%$S(D+C)\cap  C\cap C'\subset S(C)$.
%
\end{remark}
\begin{remark}
\label{easy}
Fix  $D\in \DXO$ and let $S=S(D)$.
Denote by 
$
\nu_S:\XS\to X
$
the  normalization of $X$ at the points in $S$.
Then $\XS$ is the disjoint union of the level curves of $D$, i.e.
$\XS=\coprod_{m\in \Z} D_m$.

For any regular smoothing $f:\X \to B$ of $X$ consider the twister
$D_f=\O_{\X}(D)_{|X}$ and its pull-back 
$\nu_S^*(D_f)$ 
to $\XS$. Observe that $\nu_S^*(D_f)$  does not depend on the choice of $f$;
indeed
if  $f'$  is another regular smoothing of $X$,
then $D_f$ and $D_{f'}$ have the same restrictions
to every curve $D_m$, hence their pull-backs to $\XS$ coincide.
Therefore  we can use the following non-ambiguous notation 
\begin{equation}
\label{}
\O_{\XS}(D):=\nu_S^*(D_f)=\nu_S^*\O_{\X}(D)_{|X}
\end{equation}
for any $f:\X\to B$ as above.
\end{remark}
The following  result due to Esteves and Medeiros, in \cite{EM},
 characterizes twisters among all line bundles on a nodal curve $X$.

\

\noindent
{\bf Corollary 6.9 in \cite{EM} p.297}
{ \it Fix $X$ a nodal curve,
 $D\in \DXO$ and $S=S(D)$; let $N\in \Pic X$.  
If $\nu_S^*N\cong \O_{\XS}(D)$ there exists a regular smoothing $f$ of $X$ such that
$N\cong  D_f$.}

\

\noindent
As we  mentioned, the converse also holds.
The   language
 used in \cite{EM}, section 6, is  different from ours; here is a small dictionary. 
 Our $D$ is $\tau_1C_1+\ldots +\tau_m C_m$ in 6.9 of \cite{EM}.
$\Upsilon$ in \cite{EM} is the set of irreducible components of the curve.
%(called $C$ in \cite{EM}),
%so that $\Z_{\Upsilon}$ corresponds to our $\Z^{\gamma}$.
%$\tau \in \Z^m$ in \cite{EM} corresponds to our $D\in \DXO$.
The partition $\mathcal{P}$ of $\Upsilon$ corresponds to our level expression of $D$
so that a subset $I$ of $\mathcal{P}$ corresponds to a non empty level curve $D_m$.
Condition 6.9.1 is $\nu_S^*N\cong \O_{\XS}(D)$ above.
Finally, \cite{EM} assumes  characteristic $0$, which is not needed for the
proof given to this result.
\subsubsection{A useful graph}
\label{graph}
Mantaining the hypothesis and notation above, we now
introduce the graph $\Gamma(S)$,
whose vertices are the connected components
of $\XS$ and  whose edges correspond to $S$.
An edge $e$ joins the two verteces corresponding to the two components
passing through the node represented by $e$.
So, 
$\Gamma(S)$ is the connected graph
obtained from the standard dual graph of $X$ by contracting to a point all the 
edges corresponding to nodes  not in $S$.
Let $b(S)=b_1(\Gamma(S), \Z)$ be its first Betti number,
so that 
$$
b(S)= \#S +1 -\#({\text {connected components of }} \XS ).
$$

\begin{cor}
\label{cortw} Let  $\mt \in \Lambda _X$ and 
 $S=S(\mt)$.  Then
there are  bijections 
$$
\tw^{\mt}X\leftrightarrow (\nu_S^*)^{-1}(\O_{\XS})\leftrightarrow (k^*)^{b(S)} 
$$
where $\Pic X \stackrel{\nu_S^*}{\la}\Pic \XS$ is the pull-back map.
\end{cor}
\begin{proof}
 As observed in \ref{easy}
(also in 6.9 of \cite{EM}) we have an
injection $\tw^{\mt}X\ha (\nu_S^*)^{-1}(\O_{\XS})$
 for any fixed $f$,  mapping $D_{f'}\in \tw^{\mt}X$  to $D_f\otimes
D_{f'}^{-1}$.  Surjectivity of such an injection follows from 
6.9 of \cite{EM}. 

The second  bijection (well known), follows from
the exact sequence
$$
1\la (k^*)^{b(S)} \la \Pic X \stackrel{\nu_S^*}{\la}\Pic \XS \la 1.
$$
\end{proof}

\begin{defi}
\label{sumtail}
Let $Q\subset X$ be a (connected) complete subcurve. We say that $Q$ is a {\it tail} of $X$ if
$Q\cap Q'$ is a separating node of $X$; we say that the node $Q\cap Q'$ {\it generates} the tails $Q$ and $Q'$.

Let $D\in \DXO$; we say that $D$ is a {\it sum of tails} if there is an expression
$D=\sum m_iQ_i +nX$ where the $Q_i$ are tails of $X$.

Let ${\ov D}\in \DX$, then  ${\ov D}$ is a {\it sum of tails} if any of its representative 
in $\DXO$ is.
Let $T\in \tw X$, we say that $T$ is a {\it sum of tails} if $\supp (T)$ is.
\end{defi}
\begin{remark}
\label{gp}
The set $\DXt$ of elements in $\DX $ that are sums of tails is clearly  a subgroup.
Let $\degt$ be its image via the multidegree map:
\begin{equation}
\label{degt}
\degt:=\mdeg (\DXt) \subset \Lambda _X.
\end{equation}
Thus $\degt$ is the group of multidegrees of  sums
of tails.
\end{remark}

Recall that 
 $\singsep$ is set of separating nodes of $X$. If 
$r\in \singsep$ and 
$Q$ is  one of the two tails
generated by  $r$, it is clear that $S(Q)=\{ r\}$. More generally, we have
\begin{lemma} 
\label{taillm} Let $D\in \DXO$.
$D$ is a sum of tails if and only if $S(D)\subset \singsep$.
(Equivalently: let $\mt\in\Lambda_X$.  \   $\mt \in \degt \Leftrightarrow S(\mt)\subset \singsep$)
\end{lemma}
\begin{proof} We can work modulo adding to $D$ a multiple of $X$.
We can assume that $D\neq nX$, otherwise the statement is obvious since $S(D)=\emptyset$. 

Let $D=\sum_1^l m_iQ_i$ be a sum of tails and let us use induction on $l$
to prove that
$S(D)\subset
\singsep$. If $l=1$ then $D=mQ$
with $Q$ a tail, hence $S(D)=Q\cap Q'=\{r\}$ where $r$ is the separating node generating $Q$.

If $l>1$ write $D=\sum _1^{l-1}m_iQ_i+m_lQ_l$; 
by \ref{sumD} we have $$S(D)\subset S(\sum _1^{l-1}m_iQ_i)\cup S(m_lQ_l);$$
by induction, each of the two sets  on the right lies in $\singsep$ so we are done.

Conversely, assume that $S(D)\subset \singsep$. Introduce the level expression $D=\sum_{m\in \Z} mD_m $
defined in (\ref{Ddec}); 
then every subcurve $D_m$ meets its complement $D_m'$ in separating
nodes of $X$, by (\ref{S(D)}). Obviously it suffices to prove that every connected component of $D_m$ is
a union of tails. To do that, it suffices to prove that if $Z$ is a connected
subcurve of $X$
 such that $Z\cap Z'\subset \singsep$,
then either $Z$ or $Z'$ is a union of tails. 

Use induction on $\#Z\cap Z'$.
If $\#Z\cap Z'=1$ then $Z$ and $ Z'$ are both tails. Let $\#Z\cap Z'>1$ and let
$r\in Z\cap Z'$. Since $Z$ is connected $Z$ is  contained in one of the two tails $Q_r$ and $Q_r'$
generated by $r$, say $Z\subset Q_r$. Let $Y:=\overline{Q_r\smallsetminus Z}$ so that $Z'=Y\cup
Q_r'$; then $Z\cap Y$ is made of separating nodes for $Q_r$ 
%(obviously
%the set of separating nodes of $Q_r$ equals $Q_r \cap\singsep\smallsetminus r$)
and of course $\#Z\cap Y=\#Z\cap Z'-1$. Therefore we can use
induction and conclude that either $Z$ or $Y$ is a union of tails of $Q_r$.
If $Y$ is a union of tails of $Q_r$ then $Y\cup Q_r'=Z'$ is a union of tails of $X$ and we are done.
If $Z$ is union of  tails of $Q_r$ then $Z$ is actually a tail in $Q_r$ ($Z$ is connected)
hence $Y$ is also a tail of $Q_r$ and hence (arguing as before) $Z'$ is a union of tails of $X$.
\end{proof}

\begin{cor}
\label{indtwi}  
Fix the curve $X$ and $\mt \in \Lambda_X$. Then
$$\#\tw ^{\mt}X=1 \Leftrightarrow S(\mt)\subset \singsep \Leftrightarrow \mt \in \degt .$$
(Equivalently, with the terminology of \ref{Df}: let $D\in \DXO$;
then $D_f$ does not depend on   $f$   if and only if $S(D)\subset
\singsep$ if and only if $D$ is a sum of tails.)
\end{cor}
\begin{proof}
The second double arrow is \ref{taillm}, so  we just need to prove the first.
  Fix $D\in \DXO$ such that $\mdeg D = \mt$ and let $S:=S(\mt)=S(D)$.

By \ref{cortw} we have  $\#\tw ^{\mt}X=1$ if and only if $b(S)=0$
where $b(S)$ is the first Betti number of the graph $\Gamma(S)$ defined in \ref{graph}. 
 Therefore
$b(S)=0$ if and only if $\Gamma(S)$ is a tree
if and only if $S\subset \singsep$.
\end{proof}

\section{Proof of the main Theorem}
\label{proof}
Fix a  curve $X$, then by Proposition~\ref{key} every Abel map of $X$ is of type
$\ark:\dfib^d \to \Pic^dX$,   for some $\rep$ and $f$.
To say that $\ark$ is natural (i.e. independent of the choice of $f$) is to say that for every
 regular smoothing $f'$ of $X$ we have  $\ark = \alpha^{d,\rep}_{f',X}$.
We begin with a preliminary characterization.

\begin{lemma}
\label{critind} 
The   map $\ark$ is natural if and only if 
 $\md - \rep([\md])\in \degt$
for every $\md\in \pg$.
\end{lemma}
\begin{proof}
Let us revisit  the factorization (\ref{arfact})  of 
$\ar$ by writing it$$
\ar:\dX ^d \stackrel{\mu}{\la} \frac{ \coprod_{\md \in \pg} \Pic_f^{\md}}{\sim_K} 
\stackrel{\cong}{\la}
\frac{ \coprod_{\md\in \pg} \Pic^{\rep([\md])}_f}{\sim_K}
$$
where $\mu=\mu_{E_d}$ is
the  moduli map of the  divisor $E_d\subset \dX^d\times_B \X$,
and the isomorphism is the restriction of $\iota_{\rep}\circ q_f$. 

The restriction of $\mu$ to the closed fiber $\dfib ^d$ is  fixed, i.e.
independent of $f$,
indeed if $p_1,\ldots, p_d\in \dfib$ then $\mu(p_1,\ldots, p_d)=\O_X(\sum p_i)$.
We obtain that
$\ark$ is independent of $f$ if and only if for every $\md \in \pg$ the isomorphism
$
\Pic^{\md}X \la \Pic ^{\rep([\md ])}X
$
is independent of $f$.
Such a map is 
 given by  tensor product by  that twister $T\in \tw_fX$ whose
multidegree satisfies
$\mdeg T =\md - \rep([\md ])$; set  
$
\mt :=\md - \rep([\md ])
$.
Recall that
$T$ is uniquely determined by $f$ and  
 $\mt$ (see  \ref{twlm} part (\ref{twuni})). 

We obtain that $\ark$ does not depend on the choice of $f$ if and only if $T$ is the same for
every
$f$, i.e. if and only if 
$
\tw^{\mt}X=\{ T\}
$.
We conclude by
\ref{indtwi}, which tells us that 
$
\tw^{\mt}X=\{ T\}
$ if and only if $\mt\in \degt$. 
\end{proof}
\begin{remark}
\label{uni}
Thus, for any fixed $d$ the set of natural $d$-th Abel maps is either
empty or in bijective correspondence with $(\degt)^{\#\pg}$. So,  if a
natural
$d$-th Abel map for $X$ 
exists, it is unique if and only if $X$
is free from separating nodes.
\end{remark}

\noindent
{\it Proof of Theorem~\ref{main}.}
We first prove that, if $\ce(X)>d$,
then $X$ admits a natural $d$-th Abel map, which is the  harder part.
Define a choice of representatives $\rep$ as follows. Pick $\delta\in \dcg^d$; if
$\delta $ does not contain any partitional representative,  the way
in which  we
define $\rep(\delta)$  does not matter. If instead $\delta$ contains some partitional representative,
we choose one of them, call it $\md^{\delta}$, and define $\rep(\delta)=\md^{\delta}$.
We claim that   the Abel map $\ark$ is natural, i.e. it does not depend on   $f$.
 To prove that it suffices to show that for every pair $\md, \md'\in \pg$ such that $\md\equiv \md'$
we have
\begin{equation}
\label{red1}
\md-\md'\in \degt.
\end{equation}
Indeed,  for any $\md\in \pg$, denoting by $\delta\in \dcg^d$ its class, by (\ref{red1}) we get
$
\md - \rep({\delta})=\md-\md^{\delta}\in \degt
$.
Hence, by \ref{critind},  $\ark$ does not depend on $f$.

Let $\mt:=\md-\md'$ and consider $S(\mt)$ (see (\ref{S(t)})). By \ref{taillm}, 
(\ref{red1}) is equivalent to
\begin{equation}
\label{}
S(\mt)\subset \singsep.
\end{equation}
In conclusion: we reduced ourselves to prove   the following statement.

\noindent
($\bullet$)
{\it Let $\mt \in \Lambda_X$ and assume that there exist $\md,\md' \in\pg$ such that
$\mt=\md-\md'$. Then $S(\mt)\subset \singsep$.}

We can of course assume that $\mt\neq 0$.
Let us introduce $D(\mt )=\sum_1^{\ell (\mt)} m_hZ_h(\mt)$ described in \ref{lvlm}; then 
$S(\mt)=S(D(\mt ))$.
We shall prove ($\bullet$)
using induction on $\ell(\mt)$.
First, we  simplify the notation by writing
$Z_0:=Z_0(\mt)$;  recall that $Z_0$,  defined in \ref{lvlm} part (\ref{c}),   is not empty.

Now on with the induction: assume $\ell(\mt)=1$. Then $D(\mt)=m_1Z_1(\mt)$ with $m_1>0$ and
$Z_0'=Z_1(\mt)$.

By contradiction assume that $S(\mt)\not\subset\singsep$; 
then $Z_0\cap Z_0'\not\subset\singsep$.
Therefore,
by the definition  of essential connectivity (\ref{kconn} form (\ref{ec2})) we have
$$
k_{Z_0}\geq \ce(X).
$$
On the other hand by \ref{lvlm} (\ref{degZ})  we have 
$$
|\mt_{Z_0}| \geq k_{Z_0}.
$$
Combining these two inequalities with the hypothesis $\ce(X)>d$  we get
$$
|\mt_{Z_0}|>d.
$$
Restricting the equality $\mt=\md-\md'$ to $Z_0$ and applying the above relation we obtain
$$
|\md'_{Z_0}|=|\md_{Z_0}|-|\mt_{Z_0}|<|\md_{Z_0}|-d.
$$
Now $|\md_{Z_0}|\leq d$ because $\md\in\pg$; hence
$|\md'_{Z_0}|<0$, which is in contradiction with the fact that $\md'\in \pg$.
This concludes the proof of the case $\ell(\mt )=1$.

Assume now that $\ell(\mt )\geq 2$. 
Again by contradiction suppose that $S(\mt )\not\subset \singsep$.
 If, as in the preceeding case,
$Z_0\cap Z_0'\not\subset\singsep$, we can argue exactly as before to obtain a contradiction.
So, suppose that
 $Z_0\cap Z_0'$ is made of separating nodes; then 
$
\mdeg Z_0\in \degt
$
by \ref{taillm}. Set
$$
\mdu : =\mt +\mdeg m_1Z_0;
$$
then ($\degt $ is a group)  $\mdu \in \degt$ if and only if $\mt \in \degt$.
By \ref{taillm} this is equivalent to saying that
$$
S(\mdu)\subset \singsep \Leftrightarrow S(\mt)\subset\singsep .
$$
To reach a contradiction we shall prove that $S(\mdu)\subset \singsep$.
We have 
$$
\mdu=\mdeg(m_1(Z_0 +Z_1(\mt))+m_2Z_2(\mt)+\ldots m_{\ell(\mt)}Z_{\ell(\mt)}).
$$
Therefore   lemma~\ref{lvlm}   applied to $\mdu$ gives
that $Z_0(\mdu)=Z_0 +Z_1(\mt)$ and  

$$
D(\mdu) = (m_2-m_1)Z_1(\mdu)+\ldots +(m_l-m_1)Z_{{\ell(\mt)-1}}(\mdu).
$$
So $\ell(\mdu)=\ell(\mt)-1$.
To be able to apply the induction hypothesis to $\mdu$ and 
  conclude that $S(\mdu)\subset \singsep$ we need to express
$\mdu$ as the difference of two elements in $\pg$.
To do that let
$$
\md '':=\md + \mdeg m_1Z_0=\md '+\mt +\mdeg m_1Z_0=\md ' +\mdu.
$$
We claim that $\md ''\in \pg$.
In fact, let $C$ be an irreducible component of $X$,
 then of course
$$
\md''_C=\md_C+m_1(C\cdot  Z_0) 
$$
and we need to show that $\md''_C\geq 0$. 
If $C\not\subset  Z_0$
 both summands on the right are nonnegative (recall that $\md\in \pg$) so we are ok.

If $C\subset  Z_0$, by (\ref{degZ}) in \ref{lvlm} we have that
$\mt_C\geq -m_1(C\cdot Z_0)$ hence 
$$
\md_C=\md'_C+\mt_C\geq -m_1(C\cdot Z_0)
$$
 (as $\md'\in \pg$).
Hence
$$
\md''_C=\md_C+m_1(C\cdot  Z_0) \geq -m_1(C\cdot Z_0)+m_1(C\cdot  Z_0)\geq 0.
$$

This shows that $\md''\in \pg$. Since $\mdu = \md''-\md'$ we can apply the induction hypothesis to $\mdu$
and obtain $S(\mdu)\subset \singsep$. We are done with one half of te proof.

Now we prove the converse, so assume that $X$
 admits a natural $d$-th Abel map. By \ref{key} this means that
  there exists $\rep$ such that
$\ark$ is independent of $f$. By \ref{critind},  $\forall \md \in \pg$ we have that
$\md - \rep([\md])\in \degt$ and hence, for every $\md,\md'\in \pg$
such that $\md\equiv \md'$
\begin{equation}
\label{hyp}
\md-\md'\in \degt
\end{equation}
(as $\md - \md'=\Bigr(\md - \rep([\md])\Bigl)+\Bigr(\rep([\md'])-\md'\Bigl)$
and both summands
$\Bigr(\ldots \Bigl)$ lie in $\degt$).

To prove that $\ce(X)>d$ it suffices to show that if $Z$ is a subcurve of $X$ such that $k_Z\leq d$, then
 \begin{equation}
\label{red3}
Z\cap Z'\subset\singsep
\end{equation}
(using definition (\ref{ec2}) of $\ce(X)$). So let $Z$ be such
a curve;
up to renaming the components of $X$ we have
$
Z=\cup_{h+1}^{\gamma}C_i$
 so that $1<h<\gamma$.
Thus
$
k_Z=\sum_1^h(C_i\cdot Z)
$.

By the assumption that $k_Z\leq d$ there exist $h$ integers $d_1,\ldots ,d_h$
such that $d_i\geq (C_i\cdot Z)$ for all $i=1\ldots, h$, and such that $\sum _1^hd_i=d$.
We have
\begin{displaymath}
(Z\cdot C_i) \left\{ \begin{array}{ll}
\leq d_i &\text{ if } i\leq h \  \text{ (by definition) }\\ 
\leq 0 &\text{ if } i>h \  \text{ (because } C_i\subset Z)\\
\end{array}\right.
\end{displaymath}
Set $\md:=(d_1,\ldots,d_{h},0\ldots, 0)$ so that   $\md \in \pg$
(because $d_i\geq (C_i\cdot Z)\geq 0$ for $i\leq h$).
Now define $\md':=\md -\mdeg Z$ so that
$\md'\equiv \md$. We have 
$$
\md'=\bigl(d_1-(C_1\cdot Z),\ldots,d_h-(C_h\cdot Z), -(C_{h+1}\cdot Z),
\ldots,  -(C_{\gamma}\cdot
Z)\bigr)
$$
therefore  $\md'\in \pg$.
By (\ref{hyp}) we obtain that
$\mdeg Z\in \degt$ which is to say that $S(Z)\subset \singsep$
(by Lemma~\ref{taillm}). But of course $S(Z)=Z\cap Z'$
and so we are done with the proof  of theorem~\ref{main}. \qed

\

We highlight some remarkable special cases of the theorem.
\begin{cor}
\begin{enumerate}
\item
\label{ct} Let $X$ be a   curve of compact type.
Then for all $d\geq 1$ every $d$-th Abel map of $X$ is natural.
\item
\label{nosep} Let $X$ be a nodal curve free from separating nodes and such that $\ce (X)>d$.
Then $X$ admits a unique natural $d$-th Abel map,   described as follows:
\begin{equation}
%\label{}
\begin{array}{lccr}
&\dfib^d&\la &\coprod_{\md \in \pg}\picX{\md}\\
&(p_1,\ldots,p_d) &\mapsto &\O_X(p_1+\ldots +p_d)
\end{array}
\end{equation}
\end{enumerate}
\end{cor}

\section{Abel maps and compactified Picard schemes}
\label{compare}
This final section is   to establish some connection with
   the vast literature on compactified Jacobians. To keep it to a  length
comparable with the rest of the paper,
various  important, interesting facts have been left out;
so, it may appear somewhat
 obscure to a reader who is not already acquainted with the theory of compactified Picard schemes.
On the other hand it will hopefully be useful to somebody wishing to apply 
or generalize our results to study and compactify Abel maps within a particular compatified
Picard scheme.

The generalized Jacobian of a  nodal curve $X$ fails to be projective,
unless $X$ is  of compact type. The problem of constructing a compactification for it,
with certain natural  properties, has been studied for a long time
by many authors and
diverse solutions exist (see \cite{cner} for a short overview and some
guide to the rich bibliography). Such compactifications  usually go under
the generic name of ``{\it  compactified Jacobians}" or of ``{\it  compactified Picard
schemes}". 

Here we   only consider compactifications
 compatible with the operation of
smoothing the curve (over a local one-parameter base, see property
P\ref{SN} below).
We shall  list, rather
 loosely, some   basic properties that are common to most 
such compactifications, 
in order  to relate our notion
of Abel maps to such constructions.

\subsubsection{}
\label{list}
Given a nodal curve $X$ a degree-$d$ compactified Jacobian of $X$ is a complete
scheme 
$\JXbar$ satisfying the following properties (and  others not needed here).

\begin{enumerate}[{\bf (P1)}]
\item
\label{SN}
 %{\it (Specialization and Naturality.)}
For every $f:\X \to B$ regular smoothing of $X$, 
  there exists a proper scheme  over $B$, 
$\pi:\Jfbar\to B$  whose generic fiber is $\Picgend$ and whose closed fiber
is $\JXbar$.

Denote by $\Jf\subset \Jfbar$ the smooth  locus of $\pi$ and
let $\JX$ be the closed fiber of $\Jf\to B$;  then $\JX$ is an open dense subset of $\JXbar$.

\item
\label{CN}
%{({\it Connection with the N\'eron model.})} 
%As a consequence of the N\'eron mapping property, 
There 
exists a canonical
$B$-morphism
$
n_f:\Jf\la \nf
$
which is the identity on the generic fiber. (This follows by the N\'eron mapping property)

\item
\label{uf}
 There exists a  $B$-morphism 
$
u_f:\Jf \la \Pic_f^d
$
such that $q_f\circ u_f=n_f$ and 
$u_f$ induces an isomorphism of $\Jf$ with its image.

\item
\label{uX}
The restriction  $u_X:\JX \to \Pic^dX$ of $u_f$ induces an isomorphism of
$\JX$ with a finite number of copies of the generalized Jacobian of $X$.
\end{enumerate}

\noindent
The rest of the section applies to all compactified Jacobians 
$\JXbar$ that  satisfy all the properties
in \ref{list}. 
For  cases when the map $n_f$ of (P\ref{CN})
is an isomorphism  see \cite{cner}.

Our definition~\ref{abdef} of
Abel maps does not involve compactified Jacobians, but only the  Picard scheme.
As we said in the introduction, one can 
 define Abel maps having a specific compactified Jacobian
as target  (see \cite{AK}, \cite{esteves}, \cite{CE},
\cite{coelho}). To compare such an approach to ours
we now introduce   an ad hoc
terminology, slightly awkward but useful to make distinctions.

\begin{defi}
\label{PA} 
Let $\JXbar$ be a compactified Jacobian for $X$ as in  \ref{list}.
\begin{enumerate}[(i)]
\item
\label{}
A map $\zeta:\dfib^d\la \JXbar $ is a {\it pre-Abel map} of degree $d$
if there exists a regular smoothing $f:\X\to B$ of $X$ and a map
$\zeta_f:\dX^d\la \Jfbar $
such that $(\zeta_f)_K=\alpha^d_K$ and $(\zeta_f)_k=\zeta$.
\item
\label{}
A pre-Abel map is called {\it nonsingular} if $\zeta (\dfib ^d)\subset \JX$.
%Otherwise it is called {\it singular}.
\end{enumerate}
\end{defi}

So, a nonsingular $\zeta$   maps the nonsingular locus of $X^d$ to the nonsingular locus
of $\JXbar$.
 For example the pre-Abel maps 
studied in \cite{AK} 
and  \cite{CE} 
(called simply ``Abel maps" in  such papers)
are nonsingular. This is almost immediate for the first  paper, since the curves
studied there are integral. 
The other paper 
explores the case of reducible stable curves, 
restricting to ``$d$-general curves" (cf. \cite{CE} 3.6 and 3.10)
when  $d\geq 2$, and thus getting nonsingular  Abel maps.

A nonsingular pre-Abel map determines an Abel map.
More precisely:
\begin{prop}
\label{connect}
 Let $\zeta:\dfib^d\la \JXbar $  be a nonsingular pre-Abel map;
the composition $u_X\circ \zeta:\dfib ^d \la \Pic^dX$ (notation of (P\ref{uX})) is an Abel
map.
\end{prop}
\begin{proof} The nonsingularity assumption enables us to define
$\beta:=u_X\circ \zeta$ and
$$
\beta_f:=u_f\circ \zeta_f:\dX ^d \la \Pic^d_f
$$
extending $\beta$ (where $\zeta_f$ is given by definition~\ref{PA}).
Thus to prove that $\beta $ is an Abel map   it remains to prove that condition
(\ref{sep}) of \ref{abdef} holds.

By   P\ref{uf} and  P\ref{uX} in \ref{list} there exists a finite set $S$ of multidegrees such that
$$
u_f:\Jf \stackrel{\cong}{\la} \frac{\coprod_{\md \in S}\picf{\md}}{\sim_K}\subset \picf{d}.
$$
We shall from now on identify $\Jf$ with the image of $u_f$ as indicated above.
We must prove that for every pair of multidegrees $\md$ and $\md'$ both contained in $S$ we have that
$\md \not\equiv\md'$.
By contradiction, suppose that there is a pair of distinct equivalent multidegrees $\md$ and $\md'$ in $S$;
let $D\in \DXO$ be such that $\mdeg D= \md' -\md$.

Pick now $\L\in \Pic \X$ having 
multidegree equal to $\md$ when restricted to $X$;
in order for such an $\L$ to exist we may need make a (\'etale) base change, but this will not affect the
argument. The moduli map of $\L$, 
$\mu_\L:B\to
\picf{d}$, has therefore image in $\picf{\md}$. Then we can view $\mu_\L$ as a map
$
\mu_\L:B\la \picf{\md}\subset \Jf .
$

Set $\L':=\L\otimes \O_\X(D)$, so that $\L'$ restricted to $X$ has  multidegree $\md'$. Then,
arguing  as we did for $\L$,  the moduli map of $\L'$ is 
$
\mu_{\L'}:B\la \picf{\md'}\subset \Jf .
$
Now $\mu_{\L}$ and $\mu_{\L'}$ are different maps coinciding on the generic point of $B$
(as $\L_{|\gen}=\L'_{|\gen}$). We have
 contradicted the fact that $\Jf\to B$ is separated 
(by  P\ref{SN}).
\end{proof}

\begin{example}
\label{CE}
Consider again the  Abel maps   studied in \cite{CE}.
If $d=1$ our Theorem~\ref{main} implies that
they are all natural, which was already known by the  direct proof of 
loc.cit. 4.10 and 5.13.
If   $d\geq 2$ the  question of which of them  are natural
is interesting (see  loc.cit. 3.14) and open. 
\end{example}

\begin{example}
%\label{}
 Let $X=C_1\cup C_2$ with $\#C_1\cap C_2 =\delta \geq 2$. 
Let $\JXbar$ be a compactified degree-$d$ Jacobian
having at most $\delta -1$ irreducible components
(this definitely occurs,  see 
\cite{cner} 6.5). 
Then for every $d\geq \delta$ a
pre-Abel map
$\dfib^d \to \JXbar$ is  singular.

To prove that, observe first that  $\#\dcg^d=\delta$,  so that $\delta$ is equal to the number of
connected components of the N\'eron model (i.e. of $N_X^d$).
The assumption on the number of irreducible components of $\JXbar$,  hence of $\JX$,
implies that the image of   $n_f:\Jf \to \nf$ does not intersect some connected components
of $N_X^d$. 

Now, if $d\geq \delta$ every multidegree class in $\dcg^d$
contains some partitional representative,
therefore the image of the map $\AN$ (see (\ref{Nmap})) does intersect every connected component
of 
$N^d_X$.
If there were a nonsingular pre-Abel map $\zeta:\dfib^d \la \JX$, then we could factor
$\AN=n_f\circ \zeta_f$, which is impossible.
\end{example}

\

\noindent Lucia  Caporaso  caporaso@mat.uniroma3.it\  \\\
 Dipartimento di Matematica, Universit\`a Roma Tre\\\ 
Largo S.\ L.\ Murialdo 1 \  \  00146 Roma - Italy
\end{document}